\documentclass[onefignum,onetabnum,final]{siamart171218}

\usepackage{tikz}
\usetikzlibrary{cd}
\usepackage{lipsum}
\usepackage{amsfonts}
\usepackage{graphicx}
\usepackage{epstopdf}
\usepackage{algorithmic}
\ifpdf
  \DeclareGraphicsExtensions{.eps,.pdf,.png,.jpg}
\else
  \DeclareGraphicsExtensions{.eps}
\fi

\usepackage[top=3.375cm,
            bottom=3.375cm,
            inner=4.275cm,
            outer=4.275cm]{geometry}

\newsiamremark{remark}{Remark}
\newsiamremark{hypothesis}{Hypothesis}
\crefname{hypothesis}{Hypothesis}{Hypotheses}
\newsiamthm{claim}{Claim}
\newsiamthm{assumption}{Assumption}
\newsiamremark{example}{Example}

\headers{Contour Deformation}{D. S. Winterrose}

\title{A holomorphic mapping property \\ of analytic pseudo-differential operators}

\author{David Scott Winterrose \thanks{Department of Applied Mathematics and Computer Science, Technical University of Denmark, Kgs. Lyngby, Denmark (\email{dawin@dtu.dk}).}}

\usepackage{amsopn}

\ifpdf
\hypersetup{
  pdftitle={A holomorphic mapping property of analytic pseudo-differential operators},
  pdfauthor={David Scott Winterrose}
}
\fi

\newcommand{\dbar}{d\hspace*{-0.08em}\bar{}\hspace*{0.1em}}
\newcommand{\Op}{\textnormal{Op}}

\newcommand{\supp}{\textnormal{supp}}

\newcommand{\Real}{\textnormal{Re}\,}
\newcommand{\Imag}{\textnormal{Im}\,}
\newcommand{\vol}{\textnormal{vol}}
\newcommand{\dist}{\textnormal{dist}}

\begin{document}

\maketitle

\begin{abstract}
We study the holomorphic extendibility of $\Op(p)u$, when $p$ is an analytic symbol, and explicit information is available on the domains of holomorphic extendibility of both $p$ and $u$. By a contour deformation argument, we obtain a precise local estimate of the domain of holomorphy of $\Op(p)u$ in terms of the information on $p$ and $u$.
\end{abstract}



\section{Introduction}
In \cite{Karamehmedovic2015} Karamehmedovi\'c defines a class of analytic symbols. This class forms a subspace of the analytic-type symbols in the sense of Tr\`eves \cite{TI}.
The aim was to obtain holomorphic mapping properties for the associated operators, and apply them to Calder\'on projectors in a (local) Helmholtz-type Dirichlet problem, where the boundary is a piece of a hyperplane.\\

In this way, Karamehmedovi\'c then constructs the Dirichlet-to-Neumann map, and obtains a result on how well it preserves domains of holomorphic extendibility. 
That is, how far Neumann data extends given this information about Dirichlet data, and in fact vice-versa, by the same system of equations for the Calder\'on projectors.
It was done by showing that the symbols of the Calder\'on projectors are of that class. The class of the "analytic symbols" was first introduced by Boutet de Monvel in \cite{BoutetDeMonvel1972}, and \cite{Karamehmedovic2015} essentially reuses these, but introduces constraints \cite[pp. 3-4, Definition 2.1]{Karamehmedovic2015}. 
The domains obtained in \cite[pp. 10, Theorem 2.9]{Karamehmedovic2015} are larger than those we get here, 
and \cite{Karamehmedovic2015} has the advantage of being adapted to poly-rectangular shapes.\\

The aim of this paper is to remove the strong constraints on the symbols in \cite{Karamehmedovic2015}, and reduce them to analytic symbols in the sense of Tr\`eves \cite[pp. 262, Definition 2.2]{TI}. In the process, we will also obtain a general domain-of-extension mapping theorem.
It appears in Winterrose \cite{Winterrose2021}. Let $n\in \mathbb{N}$ be the dimension throughout.

\section{Notation} 
Denote by $S^d(\mathbb{R}^n \times \mathbb{R}^n)$ order $d\in \mathbb{R}$ $(1,0)$ H\"ormander symbols. These are the $p\in C^\infty(\mathbb{R}^n \times \mathbb{R}^n)$ satisfying, for any $\alpha,\beta \in \mathbb{N}^n_0$, the estimates
\begin{align*}
\sup_{(x,\xi)\in \mathbb{R}^n \times \mathbb{R}^n} \langle \xi \rangle^{|\alpha|-d} | \partial_x^\beta \partial_\xi^\alpha p(x,\xi) | < \infty,
\end{align*}
where we use the notation $\langle \xi \rangle = (1+ |\xi|^2)^\frac{1}{2}$ for $\xi \in \mathbb{R}^n$, and put $S^{-\infty}= \cap_{d\in \mathbb{R}}S^d$.
Associated to $p$ is $\Op(p)$, 
defined via the Fourier transform $\mathcal{F}$ on $u\in C^\infty_0(\mathbb{R}^n)$ by 
\begin{align*}
\Op(p)u(x) = \frac{1}{(2\pi)^n}\int_{\mathbb{R}^n} e^{i x \cdot \xi} p(x,\xi) \mathcal{F}u(\xi) \, d\xi
\quad
\textnormal{for all}
\quad
x\in \mathbb{R}^n,
\end{align*}
which we will later write as an oscillatory integral, regularized by using Gaussians. We use the notation $\dbar\xi = (2\pi)^{-n} d\xi$ for the scaled standard Lebesgue measure $d\xi$. Finally, $B(x,r)$ denotes the open ball in $\mathbb{R}^n$ with center at $x\in \mathbb{R}^n$ and radius $r>0$, and $\mathcal{E}'(\mathbb{R}^n)$ is the space of compactly supported distributions.

\newpage
\section{Contour Deformation} 
\begin{theorem} \label{thm:ContourDeformation}
Fix $R>0$, $\epsilon>0$, and $p\in S^d(\mathbb{R}^n \times \mathbb{R}^n)$ a symbol with $d\in \mathbb{R}$. Assume $p|_{B(0,r_0)\times \mathbb{R}^n}$ extends holomorphically into $(B(0,r_0)+iB(0,\delta_0))\times W_\epsilon$, where
\begin{align*} 
W_{\epsilon} = \{ \zeta \in \mathbb{C}^n \, | \, |\Imag \zeta| < \epsilon |\Real \zeta | \} \cap \{ \zeta \in \mathbb{C}^n \, | \, |\Real \zeta | > R \},
\end{align*}
and satisfies 
\begin{align*} 
\sup_{(x,\zeta) \in K\times W_\epsilon }\langle \Real \zeta  \rangle^{-d} |p(x,\zeta)| < \infty
\quad
\textnormal{for any}
\quad 
K\subset\subset 
B(0,r_0)+iB(0,\delta_0).
\end{align*}
Let $u\in C^\infty_0(\mathbb{R}^n)$. Suppose $u|_{B(0,r)}$ extends holomorphically into $B(0,r)+iB(0,\delta)$. Choose $r>r'>0$ and $\delta \geq \delta' >0$ so that 
\begin{align*} 
\frac{\delta'}{r-r'} < \epsilon.
\end{align*}
Then $\Op(p)u|_{B(0,\min\{r',r_0\})} $ likewise extends to $B(0,\min\{r',r_0\})+iB(0,\min\{\delta',\delta_0\})$. 
\end{theorem}
This result is similar to \cite[Theorem 2.9]{Karamehmedovic2015}, but without extra constraints on $u$ and $p$. In particular, $\Op(p)u $ is real-analytic on $B(0,\min\{r',r_0\})$, as is well-known \cite{TI}.\\

A deformation of $\mathbb{R}^n \times \mathbb{R}^n$ into $\mathbb{C}^n \times \mathbb{C}^n$ allows us to continue $\Op(p)u$ explicitly. 
The main idea is to split the oscillatory integral, and apply Stokes' theorem.\\

\begin{proof} 
Take $\chi_2\in C^\infty_0(\mathbb{R}^n)$ to be $1$ on $ \overline{B(0,2R)}$ 
but $\chi_2(\xi)\in [0,1)$ for $\xi \not\in \overline{B(0,2R)}$. 
Let $\chi_1\in C^\infty_0(B(0,r))$ be a cutoff with $\chi_1(y)=1$ when $y\in B(0,r'')$, else in $[0,1)$, where $r>r''>r'$ are chosen so that
\begin{align*}
\frac{\delta'}{r-r'}<\frac{\delta'}{r''-r'} < \epsilon.
\end{align*}
Now let $\sigma : [0,1] \times \mathbb{R}^n \times \mathbb{R}^n \to \mathbb{C}^n \times \mathbb{C}^n$ be defined by
\begin{align*}
(t,y,\xi) \mapsto 
\Big(y-it\delta' \chi_1(y)(1-\chi_2(\xi)) \frac{\xi}{|\xi|}, \xi - i t \frac{\delta' (1-\chi_1(y))}{r''-r'} (1-\chi_2(\xi))|\xi| \frac{y}{|y|} \Big),
\end{align*}
and let $w$ and $\zeta$ denote the first and second  $\mathbb{C}^n$ components of this $\sigma$, respectively.
This type of $\sigma$ is used by Boutet de Monvel in \cite[pp. 243-245]{BoutetDeMonvel1972} with sparse details.
Let us put  
\begin{align*}
\mathcal{C}(t) = \sigma(\{t\} \times \mathbb{R}^n \times \mathbb{R}^n  )
\quad
\textnormal{for all}
\quad
t\in [0,1].
\end{align*}
Under the $\sigma$ deformation, if $\chi_2(\xi) = 0$ and $|\Real(x)|<r'$, we get
\begin{align*}
\Real (i(x-w)\cdot \zeta)
&=
-\xi \cdot \Big(\Imag(x) + t\delta'\chi_1(y)\frac{\xi}{|\xi|}\Big)
+
t \frac{\delta' (1-\chi_1(y))}{r''-r'} |\xi| \frac{y}{|y|}
\cdot (\Real(x) - y) \\
&\leq -|\xi|
\Big(\frac{\xi}{|\xi|} \cdot\Imag(x) + t\delta'\chi_1(y) + t\frac{\delta' (1-\chi_1(y))}{r''-r'} \Big(|y| -|\Real(x)|   \Big) \Big) \\
&\leq -|\xi|
\Big(-|\Imag(x)| + t\delta' \chi_1(y) + t\frac{\delta' (1-\chi_1(y))}{r''-r'} \Big(|y| -|\Real(x)|   \Big) \Big) \\
&\leq -|\xi|
\Big( t\delta' -|\Imag(x)| \Big).
\end{align*}
It will ensure that deformations by $\sigma(t,\cdot,\cdot)$ give convergent integrals for $|\Imag(x)|<t\delta$.
Take $x\in B(0,r')+itB(0,\delta')$, 
and fix $\rho > 2R$ and $1\geq t_2>t_1\geq 0$. Put 
\begin{align*}
Q(\rho) = (t_1,t_2) \times \mathbb{R}^n \times (B(0,\rho) \setminus \overline{B(0,2R)}).
\end{align*}
Then $\sigma$ is injective on $Q(\rho)$, because $\Real w = y$ and $\Real \zeta = \xi$ force uniqueness of $(y,\xi)$, which, by definition of $\sigma$, shows that $t$ must also be unique as long as $\xi \not\in \overline{B(0,2R)}$. Similarly, $\Real \zeta = \xi$ and $|\Imag \zeta| \leq t \epsilon |\xi|$ shows that
\begin{align*}
\sigma(\overline{Q(\rho)}) \subset \mathbb{C}^n \times W_\epsilon
\quad
\textnormal{for all}
\quad
\rho >2R.
\end{align*} 
In the following, we will put $dw = dw_1 \wedge \cdots \wedge dw_n$ and  $\dbar \zeta = (2\pi)^{-n} d\zeta_1 \wedge \cdots \wedge d\zeta_n$.
Define for $(w,\zeta)\in \mathbb{C}^n \times W_\epsilon$ the $2n$-form
\begin{align*}
\mu_x = G_x(w,\zeta) \, dw \wedge \dbar \zeta  = e^{i\zeta\cdot(x-w)} p (x,\zeta) u(w)  \, dw \wedge \dbar \zeta,
\end{align*}
where $\sigma^*\mu_x$ is smooth and compactly supported in $\overline{Q(\rho)}$, and 
\begin{align*}
d\mu_x = \sum_{j=1}^n \partial_{\overline{w}_j} \big[ e^{i\zeta\cdot(x-w)} p (x,\zeta) u(w) \big] \, d\overline{w}_j  \wedge \, dw \wedge \dbar \zeta.
\end{align*}
Then $\sigma^* d \mu_x|_{Q(\rho)}=0$, by holomorphy in $y\in B(0,r)$, and since $w \in \mathbb{R}^n$ if $y\not\in B(0,r)$.\\

Next, we show $\sigma$ is an injective immersion, and calculate its pullbacks at fixed $t$. In order to shorten expressions, we write
\begin{align*}
s(y,\xi) &= \delta'\chi_1(y)(1-\chi_2(\xi)),  \\
\eta(y,\xi) &= \frac{\delta' (1-\chi_1(y))}{r''-r'} (1-\chi_2(\xi)).
\end{align*}
Then we can calculate
\begin{align*}
dw_j &= 
dy_j-it \frac{\xi_j}{|\xi|} \sum_{k=1}^n \partial_{y_k} s(y,\xi) dy_k 
-it \sum_{k=1}^n \partial_{\xi_k} \Big( s(y,\xi) \frac{\xi_j}{|\xi|} \Big) d\xi_k  -i s(y,\xi) dt, \\
d\zeta_j &= 
d\xi_j-it \frac{y_j}{|y|}  \sum_{k=1}^n \partial_{\xi_k} \Big( \eta(y,\xi) |\xi|  \Big) d\xi_k 
 -it |\xi|\sum_{k=1}^n \partial_{y_k} \Big( \eta(y,\xi) \frac{y_j}{|y|} \Big) dy_k - i\eta(y,\xi) dt.
\end{align*}
It follows then that the real Jacobian of $\sigma$ has rank $2n$, so $\sigma$ is an injective immersion. But with $t$ kept fixed, $\det d_{(y,\xi)}\sigma(t,y,\xi)$ equals the determinant of
\begin{align*} 
\renewcommand\arraystretch{2.5}
\begin{bmatrix}
\Big[ \delta_{kj} - it \frac{\xi_j}{|\xi|} \partial_{y_k} s(y,\xi) \Big]_{k,j=1}^n & \Big[ -it \partial_{\xi_k} ( s(y,\xi) \frac{\xi_j}{|\xi|} ) \Big]_{k,j=1}^n \\
\Big[ -it |\xi| \partial_{y_k} ( \eta(y,\xi) \frac{y_j}{|y|} ) \Big]_{k,j=1}^n & \Big[\delta_{kj} -it \frac{y_j}{|y|}  \partial_{\xi_k} ( \eta(y,\xi) |\xi| ) \Big]_{k,j=1}^n
\end{bmatrix},
\end{align*}
which is bounded in $(y,\xi)$, unlike the determinant in \cite[pp. 6, Proof of Theorem 2.6]{Karamehmedovic2015}.
This term appears when pulling back 
\begin{align*}
\sigma(t,\cdot,\cdot)^*(dw \wedge \dbar\zeta) = \det d_{(y,\xi)}\sigma(t,y,\xi) \, dy \wedge \dbar \xi.
\end{align*}

\newpage
Using the above, we can now, without convergence issues, apply Stoke's theorem.
Stokes' theorem for manifolds with corners \cite[Theorem 16.25]{Lee2012} applied to $\overline{Q(\rho)}$ gives
\begin{align*}
0 = \int_{Q(\rho)} \sigma^*d\mu_x = \int_{Q(\rho)} d(\sigma^*\mu_x) = \int_{\partial Q(\rho)} \sigma^*\mu_x.
\end{align*}
Also, by the above estimate, there is some $C>0$ such that
\begin{align*}
| [(G_x\circ \sigma) \det \, d_{(y,\xi)} \sigma](t,y,\xi)|
\leq C
e^{- |\xi|(t\delta' - |\Imag(x)|)} 
\langle \xi \rangle^{d}
1_{\supp(u)}(y),
\end{align*}
which ensures existence of $\int_{\mathcal{C}(t)} \mu_x$ when $t>0$.
If $t=0$, it is meaningful if $p\in S^{-\infty}$, but $x$ must then have zero imaginary part. The aim is to show equivalence with $t=1$.
Let $\sigma_\rho : [t_1,t_2] \times \mathbb{R}^n \times \mathbb{S}^{n-1} \to\mathbb{C}^n \times \mathbb{C}^n$ be defined by
\begin{align*}
(t,y,\omega) 
\mapsto
\Big(y-it\delta' \chi_1(y) \omega, \rho \Big[ \omega - i t \frac{\delta' (1-\chi_1(y))}{r''-r'} \frac{y}{|y|} \Big] \Big).
\end{align*}
Similarly, if $x\in B(0,r')$, we get $C'>0$ such that
\begin{align*}
|[(G_x\circ \sigma_\rho) \det ( d \sigma_\rho )](t,y,\omega)|
\leq
C' e^{- \rho t\delta'} 
\langle \rho \rangle^{d+n}
1_{\supp(u)}(y),
\end{align*}
and as $\sigma(t,y,\xi)=(y,\xi)$ for $\xi\in \overline{B(0,2R)}$, $\sigma^*\mu_x$ vanishes on $(t_1,t_2) \times \mathbb{R}^n \times \partial B(0,2R)$. Combining integrals of opposite orientation, we get
\begin{align*}
\int_{\mathcal{C}(t_2)} \mu_x - \int_{\mathcal{C}(t_1)} \mu_x 
&= \lim_{\rho \to \infty}\int_{t_1}^{t_2} \int_{y\in \mathbb{R}^n} 
\int_{\xi\in\partial B(0,\rho)} (\sigma^* \mu_x)(t,y,\xi) \\
&= \lim_{\rho \to \infty}\int_{t_1}^{t_2} \int_{\mathbb{R}^n} 
\int_{ \mathbb{S}^{n-1}} 
[(G_x\circ \sigma_\rho)  \det ( d \sigma_\rho ) ] (t,y,\omega) \,
\vol_{\mathbb{S}^{n-1}}(\omega)
\, dy \, dt,
\end{align*}
where the integrand is compactly supported in $y$, bounded as above for every $\rho>R$.
It follows that the limit is zero, and we obtain that
\begin{align*}
\int_{\mathcal{C}(t_2)} \mu_x = \int_{\mathcal{C}(t_1)} \mu_x
\quad
\textnormal{if}
\quad
x\in B(0,r').
\end{align*}
Pick $t_0 \in (0,1)$ so that $\mathcal{C}(t_0)\subset \mathbb{C}^n \times W_{\frac{1}{2}}$.
By dominated convergence, we get
\begin{align*}
\Op(p)u(x) 
&=
\lim_{\lambda \to 0}
\int_{\mathbb{R}^n} \int_{\mathbb{R}^n} 
e^{i\xi\cdot(x-y)} [e^{-\lambda^2 \xi\cdot \xi}p (x,\xi)] u(y) \, dy \,\dbar\xi \\
&=
\lim_{\lambda \to 0}
\int_{\mathcal{C}({t_0})} 
e^{i\zeta\cdot(x-w)} [e^{-\lambda^2 \zeta\cdot \zeta} p(x,\zeta)] u(w)  \, dw \wedge \dbar \zeta \\
&=
\int_{\mathcal{C}({t_0})} 
e^{i\zeta\cdot(x-w)}  p(x,\zeta) u(w)  \, dw \wedge \dbar \zeta \\
&=
\int_{\mathcal{C}(1)} 
e^{i\zeta\cdot(x-w)} p (x,\zeta) u(w)  \, dw \wedge \dbar \zeta,
\end{align*}
which makes sense, because if $\lambda\in \mathbb{R}$, we have
\begin{align*}
|e^{-\lambda^2 \zeta \cdot \zeta}| \leq 
e^{-\lambda^2(|\Real \zeta |^2-|\Imag \zeta |^2)}
\leq
e^{-\frac{1}{2}\lambda^2|\Real \zeta |^2}
\quad
\textnormal{if}
\quad 
|\Imag \zeta |< \frac{1}{2}|\Real \zeta |.
\end{align*}
But now the last integral extends holomorphically in $x$ to the right open set. 
\end{proof}

\newpage
Note that for $y\not \in B(0,r)$ the function $u$ in $\mu_x$ may fail to extend holomorphically. But this is not an issue, as deformation then only happens in the $\zeta$-variable.

\begin{corollary} \label{cor:ContourDeformationDistributions}
The conclusions of Theorem~\ref{thm:ContourDeformation} hold if $u\in \mathcal{E}'(\mathbb{R}^n)$.
\end{corollary}
\begin{proof}
First pick a $\chi \in C^\infty_0(B(0,r))$ 
such that $\chi(y)=1$ for every $y\in \supp(\chi_1)$. 
Define $\sigma_y : [0,1] \times \mathbb{R}^n \to \mathbb{C}^n$ by
\begin{align*}
(t,\xi) \mapsto \zeta
=
 \xi - i t \frac{\delta' (1-\chi_1(y))}{r''-r'} (1-\chi_2(\xi))|\xi| \frac{y}{|y|},
\end{align*}
and put 
\begin{align*}
\mathcal{C}_y (1) = \sigma_y(\{1\} \times \mathbb{R}^n).
\end{align*}
As before, if $\chi_2(\xi)=0$ and $|\Real(x)| < r'$, we get
\begin{align*}
\Real (i(x-y)\cdot \zeta)
&=
-\xi \cdot \Imag(x) 
+
t \frac{\delta' (1-\chi_1(y))}{r''-r'} |\xi| \frac{y}{|y|}
\cdot (\Real(x) - y) \\
&\leq -|\xi|
\Big(\frac{\xi}{|\xi|} \cdot\Imag(x) + t\frac{\delta' (1-\chi_1(y))}{r''-r'} \Big(|y| -|\Real(x)|   \Big) \Big) \\
&\leq -|\xi|
\Big( t \delta' (1-\chi_1(y))  - |\Imag(x)|  \Big).
\end{align*}
Taking $\varphi \in C^\infty_0(B(0,r'))$, we have
\begin{align*}
\langle \Op(p)u , \varphi \rangle
&=
\langle \Op(p)(\chi u) , \varphi \rangle 
+
\langle \Op(p)((1-\chi)u) , \varphi \rangle \\
&=
\langle \Op(p)(\chi u) , \varphi \rangle 
+
\int_{\mathbb{R}^n} \Big\langle u(y), K(x,y) \Big\rangle  \, \varphi(x) \, dx \\
&=
\Big\langle \Op(p)(\chi u)(x) + \langle u(y), K(x,y) \rangle, \varphi(x) \Big\rangle,
\end{align*}
where $K: B(0,r')\times \mathbb{R}^n \to \mathbb{C}$ is the smooth
kernel of $\Op(p)(1-\chi)$ on $B(0,r')$ only, and we use brackets $\langle \cdot, \cdot\rangle$ to denote the pairing of a distribution and a test function.
The action of $\Op(p)$ is understood in the distributional sense via the formal adjoint.
By \cref{thm:ContourDeformation}, $\Op(p)(\chi u)$ extends holomorphically to the tube 
\begin{align*}
T = B(0,\min\{r',r_0\}) + iB(0,\min\{\delta',\delta_0\}),
\end{align*}
provided that
\begin{align*}
\frac{\delta'}{r-r'}<\frac{\delta'}{r''-r'}<\epsilon.
\end{align*}
It follows then (see e.g. \cite[pp. 53-54, Exercise 3.14]{GG}) that $x \mapsto \langle u, K(x,\cdot) \rangle$ is smooth, and all derivatives go through the brackets, because $K$ is smooth and $u \in \mathcal{E}'(\mathbb{R}^n)$.
The same is true if $K$ extends holomorphically in $x$ to a smooth $K: T \times \mathbb{R}^n \to \mathbb{C}$. We can then simply take complex derivatives through the brackets 
\begin{align*}
\partial_{\overline{z}} \langle u, K(z,\cdot) \rangle 
= \langle u, \partial_{\overline{z}}K(z,\cdot) \rangle = 0,
\end{align*}
and $\Op(p)u$ then extends (strongly) to the holomorphic function
\begin{align*}
\Op(p)u(z) 
= \Op(p)(\chi u)(z) + \langle u, K(z,\cdot) \rangle.
\end{align*}

\newpage
It remains only to show the holomorphic extension as outlined above for $K(x,y)$.
Pick $t_0 \in (0,1)$ so that $\mathcal{C}_y (t_0)\subset  W_{\frac{1}{2}}$. We deform from $0$ to $t_0$ with $\mathcal{C}_y (t_0)\subset  W_{\frac{1}{2}}$, where $p$ is multiplied by a Gaussian symbol, and finally from $t=t_0$ to $t=1$ directly. 
An argument using Stokes' theorem shows that $K$ has the form
\begin{align*}
K(x,y)
&=
\lim_{\lambda \to 0}\int_{\mathbb{R}^n} e^{i(x-y)\cdot \xi} [e^{-\lambda^2 \xi \cdot \xi}(1-\chi)(y)p(x,\xi) ] \, \dbar \xi \\
&=
\lim_{\lambda \to 0}
\int_{\mathcal{C}_y({t_0})} 
e^{i(x-y)\cdot \zeta} [e^{-\lambda^2 \zeta \cdot \zeta}(1-\chi)(y)p(x,\zeta) ] \, \dbar \zeta \\
&=
\int_{\mathcal{C}_y(1)} 
e^{i(x-y)\cdot \zeta} (1-\chi)(y)p(x,\zeta)  \, \dbar \zeta,
\end{align*}
and $K$ vanishes unless $y\not \in \supp(\chi_1)$, in which case
\begin{align*}
\Real (i(x-y)\cdot \zeta)
\leq -|\xi| ( t \delta'  - |\Imag(x)| ).
\end{align*}
This means that the last deformed integral is absolutely convergent if $|\Imag (x)|<\delta'$, and thus $K$ extends holomorphically in $x$ to a smooth function on $T\times \mathbb{R}^n$.
\end{proof}

\begin{theorem}
Let $U\subset \mathbb{R}^n$ be open, and $U_\mathbb{C}$ be a tube-domain about $U$ in $\mathbb{C}^n$. (This means $z \in U_\mathbb{C}$ implies $\Real z \in U$ and $\Real(z)+iy \in U_\mathbb{C}$ for all $|y|\leq |\Imag(z)|$.)
Assume $p|_{U \times \mathbb{R}^n}$ extends into
$U_\mathbb{C} \times W_\epsilon$, with $W_\epsilon$ as in \cref{thm:ContourDeformation}, and
\begin{align*} 
\sup_{(x,\zeta) \in K\times W_\epsilon }\langle \Real \zeta  \rangle^{-d} |p(x,\zeta)| < \infty
\quad
\textnormal{for any}
\quad 
K\subset\subset 
U_\mathbb{C}.
\end{align*}
Let $u\in \mathcal{E}'(\mathbb{R}^n)$ be real-analytic on $U$, with $u|_U$ extending holomorphically into $U_\mathbb{C}$.
Then $\Op(p)u|_U$ extends holomorphically into 
\begin{align*}
 \{ z\in U_{\mathbb{C}} \, | \,  |\Imag z | < \epsilon \, \dist(\Real z, \partial U) \},
\end{align*}
and is independent of $R>0$ in $W_\epsilon$.
\end{theorem}
\begin{proof}
\cref{cor:ContourDeformationDistributions} is valid over any $x\in U$ by translation of $x$ to the origin.
This gives a holomorphic extension of $\Op(p)u|_{B(x,r')}$ into $B(x,r')+iB(0,\delta')$ with
\begin{align*} 
\delta' < \epsilon (\dist(x, \partial U) - r'),
\end{align*}
and by making $r'$ small, we can make $\delta'$ arbitrarily close to $\epsilon \, \dist(x, \partial U)$.
\end{proof}

\section{Remarks} 
This removes the topology needed in \cite[pp. 3-4, Definition 2.1]{Karamehmedovic2015}. It reduces the situation to symbols defined by Boutet de Monvel \cite{BoutetDeMonvel1972} and Tr\`eves \cite{TI}. 
However, the approach to the original question raised in \cite{Karamehmedovic2015} has since changed a lot, and in \cite{Karamehmedovic2021}, we will approach it via precise local convergence radius estimates instead. It should be noted that those estimates do not subsume the result in this paper. \\

The reason is that it is hard to get parametrix symbols in the same analytic class.
It works in \cite{Karamehmedovic2015} because the geometry is simple - the boundary is a piece of a hyperplane.
To overcome this, the analytic symbols are replaced with pseudo-analytic amplitudes, which have weaker conditions imposed on them - analyticity is replaced by an estimate,
and gives a way to build pseudo-analytic parametrices from formal asymptotic sums.
This can be exploited to obtain controlled convergence radius estimates.

\newpage
\section{Acknowledgements}
The author wishes to thank the anonymous reviewer whose critique and suggestions have greatly improved the paper.

\bibliographystyle{siamplain}
\bibliography{references}

\begin{thebibliography}{1}

\bibitem{BoutetDeMonvel1972}
{\sc L.~Boutet~de Monvel}, {\em {Opérateurs pseudo-différentiels analytiques
  et opérateurs d'ordre infini}}, Annales de l'institut Fourier, 22 (1972),
  pp.~229--268.

\bibitem{GG}
{\sc G.~Grubb}, {\em Distributions and Operators}, vol.~252, Springer, 2009.

\bibitem{Karamehmedovic2015}
{\sc M.~Karamehmedovi\'c}, {\em On analytic continuability of the missing
  cauchy datum for helmholtz boundary problems}, Proceedings of the American
  Mathematical Society, 143 (2015), pp.~1515--1530.

\bibitem{Karamehmedovic2021}
{\sc M.~Karamehmedovi\'c and D.~S. Winterrose}, {\em Convergence radius
  estimates for a calculus of analytic pseudodifferential operators}, to
  appear.

\bibitem{Lee2012}
{\sc J.~M. Lee}, {\em Introduction to Smooth Manifolds}, vol.~218 of Graduate
  Texts in Mathematics, Springer-Verlag, 2~ed., 2012.

\bibitem{TI}
{\sc J.~F. Trèves}, {\em {Introduction to Pseudodifferential and Fourier
  Integral Operators: Pseudodifferential Operators}}, Springer, 1980.

\bibitem{Winterrose2021}
{\sc D.~S. Winterrose}, {\em {Complexifications, Pseudodifferential Operators,
  and the Poisson Transform}}, PhD thesis, Technical University of Denmark,
  2021.

\end{thebibliography}
\end{document}